\documentclass[11pt]{article}

\begin{document}

%%-I include macros files here-%%

%%dung le's macros

\newcommand{\newc}{\newcommand}

%%index

\renewcommand{\theequation}{\thesection.\arabic{equation}}
\newc{\eqnoset}{\setcounter{equation}{0}}
\newcommand{\myref}[2]{#1~\ref{#2}}

\newcommand{\mref}[1]{(\ref{#1})}
\newcommand{\reflemm}[1]{Lemma~\ref{#1}}
\newcommand{\refrem}[1]{Remark~\ref{#1}}
\newcommand{\reftheo}[1]{Theorem~\ref{#1}}
\newcommand{\refdef}[1]{Definition~\ref{#1}}
\newcommand{\refcoro}[1]{Corollary~\ref{#1}}
\newcommand{\refprop}[1]{Proposition~\ref{#1}}
\newcommand{\refsec}[1]{Section~\ref{#1}}
\newcommand{\refchap}[1]{Chapter~\ref{#1}}

%%environments
\newcommand{\beq}{\begin{equation}}
\newcommand{\eeq}{\end{equation}}
\newcommand{\beqno}[1]{\begin{equation}\label{#1}}

\newcommand{\barr}{\begin{array}}
\newcommand{\earr}{\end{array}}

\newc{\bearr}{\begin{eqnarray*}}
\newc{\eearr}{\end{eqnarray*}}

\newc{\bearrno}[1]{\begin{eqnarray}\label{#1}}
\newc{\eearrno}{\end{eqnarray}}

\newc{\non}{\nonumber}
\newc{\nol}{\nonumber\nl}

\newcommand{\bdes}{\begin{description}}
\newcommand{\edes}{\end{description}}
\newc{\benu}{\begin{enumerate}}
\newc{\eenu}{\end{enumerate}}
\newc{\btab}{\begin{tabular}}
\newc{\etab}{\end{tabular}}

%%\newtheorem{theorem}{Theorem}
%%\newtheorem{defi}[theorem]{Definition}
%%\newtheorem{lemma}{Lemma}[section]
%%\newtheorem{rem}[lemma]{Remark}
%%\newtheorem{exam}[lemma]{Example}
%%\newtheorem{propo}[theorem]{Proposition}
%%\newtheorem{corol}[theorem]{Corollary}

%%Unmark these for uniform indexing

\newtheorem{theorem}{Theorem}[section]
\newtheorem{defi}[theorem]{Definition}
\newtheorem{lemma}[theorem]{Lemma}
\newtheorem{rem}[theorem]{Remark}
\newtheorem{exam}[theorem]{Example}
\newtheorem{propo}[theorem]{Proposition}
\newtheorem{corol}[theorem]{Corollary}

\renewcommand{\thelemma}{\thesection.\arabic{lemma}}

\newcommand{\btheo}[1]{\begin{theorem}\label{#1}}
\newc{\brem}[1]{\begin{rem}\label{#1}\em}
\newc{\bexam}[1]{\begin{exam}\label{#1}\em}
\newc{\bdefi}[1]{\begin{defi}\label{#1}}
\newcommand{\blemm}[1]{\begin{lemma}\label{#1}}
\newcommand{\bprop}[1]{\begin{propo}\label{#1}}
\newcommand{\bcoro}[1]{\begin{corol}\label{#1}}
\newcommand{\etheo}{\end{theorem}}
\newcommand{\elemm}{\end{lemma}}
\newcommand{\eprop}{\end{propo}}
\newcommand{\ecoro}{\end{corol}}
\newc{\erem}{\end{rem}}
\newc{\eexam}{\end{exam}}
\newc{\edefi}{\end{defi}}

\newc{\rmk}[1]{{\bf REMARK #1: }}
\newc{\DN}[1]{{\bf DEFINITION #1: }}

\newcommand{\bproof}{{\bf Proof:~~}}
\newc{\eproof}{{\vrule height8pt width5pt depth0pt}\vspace{3mm}}

%%symbols

\newcommand{\rarr}{\rightarrow}
\newcommand{\Rarr}{\Rightarrow}
\newcommand{\tru}{\backslash}
\newc{\bfrac}[2]{\dspl{\frac{#1}{#2}}}

%%space

\newc{\nl}{\vspace{2mm}\\}
\newc{\nid}{\noindent}

%%operations

\newcommand{\oneon}[1]{\frac{1}{#1}}
\newcommand{\dspl}{\displaystyle}
\newc{\grad}{\nabla}
\newc{\Div}{\mbox{div}}
\newc{\pdt}[1]{\dspl{\frac{\partial{#1}}{\partial t}}}
\newc{\pdn}[1]{\dspl{\frac{\partial{#1}}{\partial \nu}}}
\newc{\pdNi}[1]{\dspl{\frac{\partial{#1}}{\partial \mathcal{N}_i}}}
\newc{\pD}[2]{\dspl{\frac{\partial{#1}}{\partial #2}}}
\newc{\dt}{\dspl{\frac{d}{dt}}}
\newc{\bdry}[1]{\mbox{$\partial #1$}}
\newc{\sgn}{\mbox{sign}}

\newc{\Hess}[1]{\frac{\partial^2 #1}{\pdh z_i \pdh z_j}}
\newc{\hess}[1]{\partial^2 #1/\pdh z_i \pdh z_j}

%%function spaces

\newcommand{\Coone}[1]{\mbox{$C^{1}_{0}(#1)$}}
\newcommand{\lspac}[2]{\mbox{$L^{#1}(#2)$}}
\newc{\hspac}[2]{\mbox{$C^{0,#1}(#2)$}}
\newc{\Hspac}[2]{\mbox{$C^{1,#1}(#2)$}}
\newc{\Hosp}{\mbox{$H^{1}_{0}$}}
\newcommand{\Lsp}[1]{\mbox{$L^{#1}(\Og)$}}
\newc{\hsp}{\Hosp(\Og)}

%%greeks

\newc{\ag}{\alpha}
\newc{\bg}{\beta}
\newc{\cg}{\gamma}\newc{\Cg}{\Gamma}
\newc{\dg}{\delta}\newc{\Dg}{\Delta}
\newc{\eg}{\varepsilon}
\newc{\zg}{\zeta}
\newc{\thg}{\theta}
\newc{\llg}{\lambda}\newc{\LLg}{\Lambda}
\newc{\kg}{\kappa}
\newc{\rg}{\rho}
\newc{\sg}{\sigma}\newc{\Sg}{\Sigma}
\newc{\tg}{\tau}
\newc{\fg}{\phi}\newc{\Fg}{\Phi}
\newc{\vfg}{\varphi}
\newc{\og}{\omega}\newc{\Og}{\Omega}
%\newc{\ng}{\eta}
\newc{\pdh}{\partial}

%%Integration and Sum

\newc{\ii}[1]{\int_{#1}}
\newc{\iidx}[2]{{\dspl\int_{#1}~#2~dx}}
\newc{\bii}[1]{{\dspl \ii{#1} }}
\newc{\biii}[2]{{\dspl \iii{#1}{#2} }}
\newc{\su}[2]{\sum_{#1}^{#2}}
\newc{\bsu}[2]{{\dspl \su{#1}{#2} }}

\newcommand{\iiomdx}[1]{{\dspl\int_{\Og}~ #1 ~dx}}
\newc{\biiom}[1]{{\dspl\int_{\bdrom}~ #1 ~d\sg}}
\newc{\io}[1]{{\dspl\int_{\Og}~ #1 ~dx}}
\newc{\bio}[1]{{\dspl\int_{\bdrom}~ #1 ~d\sg}}
\newc{\bsir}{\bsu{i=1}{r}}
\newc{\bsim}{\bsu{i=1}{m}}

\newc{\iibr}[2]{\iidx{\bprw{#1}}{#2}}
\newc{\Intbr}[1]{\iibr{R}{#1}}
\newc{\intbr}[1]{\iibr{\rg}{#1}}
\newc{\intt}[3]{\int_{#1}^{#2}\int_\Og~#3~dxdt}
%%\newc{\itQ}[2]{\dspl{\int\hspace{-2.5mm}\int_{#1}~#2~dxdt}}
%%\newc{\mitQ}[2]{\dspl{\rule[1mm]{4mm}{.3mm}\hspace{-5.3mm}\int\hspace{-2.5mm}\int_{#1}~#2~dxdt}}

\newc{\itQ}[2]{\dspl{\int\hspace{-2.5mm}\int_{#1}~#2~dz}}
\newc{\mitQ}[2]{\dspl{\rule[1mm]{4mm}{.3mm}\hspace{-5.3mm}\int\hspace{-2.5mm}\int_{#1}~#2~dz}}
\newc{\mitQQ}[3]{\dspl{\rule[1mm]{4mm}{.3mm}\hspace{-5.3mm}\int\hspace{-2.5mm}\int_{#1}~#2~#3}}

\newc{\mitx}[2]{\dspl{\rule[1mm]{3mm}{.3mm}\hspace{-4mm}\int_{#1}~#2~dx}}

\newc{\mitQq}[2]{\dspl{\rule[1mm]{4mm}{.3mm}\hspace{-5.3mm}\int\hspace{-2.5mm}\int_{#1}~#2~d\bar{z}}}
\newc{\itQq}[2]{\dspl{\int\hspace{-2.5mm}\int_{#1}~#2~d\bar{z}}}

\newc{\pder}[2]{\dspl{\frac{\partial #1}{\partial #2}}}

%%variables

\newc{\ui}{u_{i}}
\newcommand{\upl}{u^{+}}
\newcommand{\umn}{u^{-}}
\newcommand{\un}{\{ u_{n}\}}

\newcommand{\uo}{u_{0}}
\newc{\voi}{v_{i}^{0}}
\newc{\uoi}{u_{i}^{0}}
\newc{\vu}{\vec{u}}

\newc{\xo}{x_{0}}
\newc{\Br}{B_{R}}
\newc{\Bro}{\Br (\xo)}
\newc{\bdrom}{\bdry{\Og}}
\newc{\ogr}[1]{\Og_{#1}}
\newc{\Bxo}{B_{x_0}}

%%Additional macros
\newc{\inP}[2]{\|#1(\bullet,t)\|_#2\in\cP}
\newc{\cO}{{\mathcal O}}
\newc{\inO}[2]{\|#1(\bullet,t)\|_#2\in\cO}

\newc{\newl}{\\ &&}

%\newc{\bilhom}{\mbox{Bil}(\mbox{Hom}(\RR^n,\RR^N))}
\newc{\bilhom}{\mbox{Bil}(\mbox{Hom}(\RR^{nm},\RR^{nm}))}
\newc{\VV}[1]{{V(Q_{#1})}}

\newc{\ccA}{{\mathcal A}}
\newc{\ccB}{{\mathcal B}}
\newc{\ccC}{{\mathcal C}}
\newc{\ccD}{{\mathcal D}}
\newc{\ccE}{{\mathcal E}}
\newc{\ccH}{\mathcal{H}}
\newc{\ccF}{\mathcal{F}}
\newc{\ccI}{{\mathcal I}}
\newc{\ccJ}{{\mathcal J}}
\newc{\ccP}{{\mathcal P}}
\newc{\ccQ}{{\mathcal Q}}
\newc{\ccR}{{\mathcal R}}
\newc{\ccS}{{\mathcal S}}
\newc{\ccT}{{\mathcal T}}
\newc{\ccX}{{\mathcal X}}
\newc{\ccY}{{\mathcal Y}}
\newc{\ccZ}{{\mathcal Z}}

\newc{\bb}[1]{{\mathbf #1}}
\newc{\bbA}{{\mathbf A}}
%%Bracket
\newc{\myprod}[1]{\langle #1 \rangle}
\newc{\mypar}[1]{\left( #1 \right)}

%%Norms

\newc{\lspn}[2]{\mbox{$\| #1\|_{\Lsp{#2}}$}}
\newc{\Lpn}[2]{\mbox{$\| #1\|_{#2}$}}
\newc{\Hn}[1]{\mbox{$\| #1\|_{H^1(\Og)}$}}

%%Misc

\newc{\cyl}[1]{\og\times \{#1\}}
\newc{\cyll}{\og\times[0,1]}
\newc{\vx}[1]{v\cdot #1}
\newc{\vtx}[1]{v(t,x)\cdot #1}
\newc{\vn}{\vx{n}}

\newcommand{\RR}{{\rm I\kern -1.6pt{\rm R}}}

%\numberwithin{equation}{section}
%\newcommand{\eproof}{\vrule height5pt width3pt depth0pt}
%\newtheorem{exam}[lem]{Example}

\newenvironment{proof}{\noindent\textbf{Proof.}\ }
{\nopagebreak\hbox{ }\hfill$\Box$\bigskip}

%%Additional macros

\newc{\itQQ}[2]{\dspl{\int_{#1}#2\,dz}}
\newc{\mmitQQ}[2]{\dspl{\rule[1mm]{4mm}{.3mm}\hspace{-4.3mm}\int_{#1}~#2~dz}}
\newc{\MmitQQ}[2]{\dspl{\rule[1mm]{4mm}{.3mm}\hspace{-4.3mm}\int_{#1}~#2~d\mu}}

\newc{\MUmitQQ}[3]{\dspl{\rule[1mm]{4mm}{.3mm}\hspace{-4.3mm}\int_{#1}~#2~d#3}}
\newc{\MUitQQ}[3]{\dspl{\int_{#1}~#2~d#3}}

\vspace*{-.8in}
\begin{center} {\LARGE\em  Global Existence and Regularity Results for Large Cross Diffusion Models  on Planar Domains.}

 \end{center}

\vspace{.1in}

\begin{center}

{\sc Dung Le}{\footnote {Department of Mathematics, University of
Texas at San
Antonio, One UTSA Circle, San Antonio, TX 78249. {\tt Email: Dung.Le@utsa.edu}
}}

\end{center}

\begin{abstract}
The global existence of classical solutions to cross diffusion systems of more than 2 equations given on a planar domain is established. The results can apply to generalized Shigesada-Kawasaki-Teramoto (SKT) and food pyramid models whose diffusion and reaction can have polynomial growth of any order.
\end{abstract}

\vspace{.2in}

\section{Introduction}\label{introsec}\eqnoset

 We consider in this paper the following system
\beqno{e1}\left\{\barr{ll} u_t=\Div(A(u)Du)+f(u,Du)& (x,t)\in Q=\Og\times(0,T),\\u(x,0)=U_0(x)& x\in\Og\\\mbox{$u=0$ on $\partial \Og\times(0,T)$}. &\earr\right.\eeq

Here, $\Og$ is a bounded domain with smooth boundary $\partial \Og$ in $\RR^2$; $u:\Og\to\RR^m$ and $f:\RR^m\times\RR^{2m}\to\RR^m$ are vector valued functions. $A(u)$ is a full matrix $m\times m$. Thus, the above is a system of $m$ equations.

The system \mref{e1} arises in many mathematical biology and ecology applications. In the last few decades, papers concerning such strongly coupled parabolic systems usually assumed that the solutions under consideration were bounded, a very hard property to check as maximum principles had been unavailable for systems in general. Most of global existence results for cross diffusion systems relied on the following local existence result of Amann.

\btheo{Amthm} (\cite{Am1,Am2}) Suppose $\Og\subset\RR^n$, $n \ge2$, with $\partial\Og$
 being smooth. Assume that \mref{e1} is normally elliptic. Let $p_0 \in (n,\infty)$ and $U_0$ be  in $W^{1,p_0}(\Og
)$. Then there exists a maximal time $T_0\in(0,\infty]$ such that the system
\mref{e1} has a unique classical solution in $(0, T_0)$ with
$$u\in C([0,T_0),W^{1,p_0}(\Og)) \cap C^{1,2}((0, T_0)\times\bar{\Og})$$
Moreover, if $T_0<\infty$ then \beqno{blowup}\lim_{t\to T^-_0}\|u(\cdot,t)\|_{W^{1,p_0}(\Og)}=\infty.\eeq\etheo

Equivalently, the classical solution $u$ will exist globally if its $W^{1,p_0}(\Og)$ norm does not blow up in finite time. This requires the existence of a continuous function $\ccC$ on $(0,\infty)$ such that \beqno{blowupz} \|u(\cdot,t)\|_{W^{1,p_0}(\Og)} \le \ccC(t) \quad \forall t\in(0,T_0).\eeq

We refer the readers to \cite{Am1} for the definition of normal ellipticity (roughly speaking, it means that the real parts of the eigenvalues of $A(u)$ are positive). The checking of \mref{blowup} is the most difficult one as known techniques for the regularity of solutions to scalar equations could not be extended to systems and counterexamples were available.  

In our recent work \cite{leGN}, see also \cite{letrans,leans}, we considered \mref{e1} on a domain in $\RR^n$ ($n\ge2$) and were able to relax the condition \mref{blowupz} by \beqno{blowupzz}\|u(\cdot,t)\|_{W^{1,n}(\Og)}\le \ccC(t) \quad \forall t\in(0,T_0).\eeq

Obviously, \mref{blowupzz} does not imply that $|u|$ is bounded so that \mref{e1} is not regularly elliptic, i.e. eigenvalues of $A(u)$ can still be unbounded. In \cite{leGN}, we only assume that \mref{e1} is uniformly elliptic, meaning the the eigenvalues of $A(u)$ are comparable.

In fact, we prove in \cite{leGN} that the global existence result in Amann's theorems holds with a much weaker version of \mref{blowupz} and \mref{blowupzz}. Namely,  one needs only to control the BMO norm of $u$ and shows that $u$ is VMO. The proof makes use several techniques from Harmonic Analysis and a generalized weighted Gagliardo-Nirenberg inequality involving BMO norms. Some mild structural conditions on the cross diffusion matrix $A(u)$ of \mref{e1} are imposed and easily verified, for examples, if $A(u)$ has a polynomial growth of order less than 5.

In this work, for planar domains ($n=2$) we will present a new and simpler proof of the described result in \cite{leGN}. Tools from Harmonic Analysis will not be needed here. More importantly, the structural conditions on $A(u)$ will be much weaker. As an example, we will allow $A(u)$ has  a polynomial growth of {\em any} order. The new result is given in \reftheo{thm0}.

The checking of \mref{blowupzz} for general $n$ is by no means an easy task. Here, for $n=2$ we will give two examples in applications where this can be done under very general assumptions.

In paticular, we will consider a class of generalized Shigesada-Kawasaki-Teramoto (SKT) models consisting of more than 2 equations. Namely, we will establish the global existence of classical solutions to the following system
\beqno{SKTgen} u_t = \Delta(\ccP(u)) + f(u,Du), \eeq where $\ccP(u),f(u,Du)$ are vector valued functions. The above system is a special case of \mref{e1} with $A(u)=\frac{\partial \ccP(u)}{\partial u}$. The (SKT) models, consisting of only 2 equations, were introduced in \cite{SKT}  using quadratic growth $\ccP$'s and Lotka-Volterra reaction $f(u)$. The global existence for the (SKT) model for 2 species on planar domains was studied in \cite{yag}. In this work, we will consider \mref{SKTgen} with $\ccP(u),f(u,Du)$ have polynomial growth of any order in $u$. More importantly, the number of species/equations can be arbitrary.

The second example is the food pyramid model which assumes that the first $k$ equations of \mref{e1} do not involve with the $j$-unknown $u_j$ if $j>k$. This model has been studied in literature (e.g., \cite{Al}) under the assumption that $A(u)$ is a constant diagonal matrix, with no cross diffusion.

\section{Preliminaries and Main Results}\eqnoset\label{res}

Throughout this paper $\Og$ is a bounded domain with smooth boundary in $\RR^2$. 
For any smooth (vector valued) function $u$ defined on $\Og\times(0,T)$, $T>0$, its temporal and spatial derivatives are denoted by $u_t,Du$ respectively. If $A$ is a $C^1$ function in $u$ then we also abbreviate $\frac{\partial A}{\partial u}$ by $A_u$.

As usual, $W^{1,p}(\Og,\RR^m)$, $p\ge1$, will denote the standard Sobolev spaces whose elements are vector valued functions $u\,:\,\Og\to \RR^m$ with finite norm $$\|u\|_{W^{1,p}(\Og,\RR^m)} = \|u\|_{L^p(\Og)} + \|Du\|_{L^p(\Og)}.$$ 

We assume the following structural conditions.

\bdes

\item[A)] $A(u)$ is $C^1$ in $u$. Moreover, there are positive constants $\llg_0,C$ and a scalar $C^1$ function $\llg(u)$  such that $\llg(u)\ge \llg_0$ for all $u\in\RR^m$. Furthermore, for any $\zeta\in\RR^{nm}$ 
\beqno{A1} \llg(u)|\zeta|^2 \le \myprod{A(u)\zeta,\zeta} \mbox{ and } |A(u)|\le C\llg(u).\eeq 

We also assume $|A_u|\le C|\llg_u|$ and   \beqno{Fghyp} |\llg_u(u)| \le C\llg(u).\eeq  

\edes

Our first main theorem under this general assumption weakens the condition \mref{blowup} of Amann's result in \reftheo{Amthm} by the condition \mref{blowup0} where we use much weaker norms.

\btheo{thm0} Suppose $\Og\subset\RR^2$ with $\partial\Og$
 being smooth. Assume A) and $U_0$ be  in $W^{1,p_0}(\Og
)$ for some $p_0>2$. Let $T_0\in(0,\infty]$ be the maximal existence time for a unique classical solution 
$u\in  C^{1,2}((0, T_0)\times\bar{\Og})$ of 
\beqno{e1z}\left\{\barr{ll} u_t=\Div(A(u)Du)+\hat{f}(u,Du)& (x,t)\in Q=\Og\times(0,T),\\u(x,0)=U_0(x)& x\in\Og\\\mbox{Boundary conditions for $u$ on $\partial \Og\times(0,T)$}. &\earr\right.\eeq

 For such $u$, asume that there are a constant $C$ and a $C^1$ function $f(u)$ such that $|f_u(u)|\le C\llg(u)$ and  $\hat{f}(u,Du)$ satisfies the following growth conditions

\beqno{FUDU1a}|\hat{f}(u,Du)| \le C\llg^\frac12(u)|Du| + f(u),\eeq 
Then if there exists a continuous function $g$ on $(0,\infty)$ such that \beqno{blowup0}\iidx{\Og}{\llg(u)|Du(x,t)|^2}\le g(t) \quad \forall t\in(0,T_0),\eeq then $T_0=\infty$. \etheo

The proof of this theorem is fairly technical and we will postpone it to the end of the paper (see \refsec{secthm0}). Of course, the assumption \mref{blowup0} is a bit stronger than \mref{blowup} but the structural conditions on \mref{e1} are much weaker than those in \cite{leGN}. In fact, we will show that the crucial hypothesis \mref{blowup0} for global existence  can be verified in many applications including the general (SKT) and pyramid systems.

For the sake of simplicity, we will consider first the case where the reaction term in \mref{e1} does not depend on $Du$. We then assume the following.

\bdes 
\item[F)] We assume that there are positive constants $\eg_0,C$ and  nonnegative $C^1$ functions $P,F:\RR^m\to \RR^+$ satisfying $F(0)=P(0)=0$ and \beqno{dfuu}|F_u(u)|\le C\llg^\frac12(u), \eeq
\beqno{Pu}|P_u(u)|\le C\llg(u)\eeq for all $u\in \RR^m$ such that
\beqno{fuu}|f(u)||u|\le \eg_0 F^2(u)+C,\eeq \beqno{fuuu}  \frac{\llg^\frac12(u)|f(u)|}{P(u)+1}\le C(F(u)+1).\eeq

\edes

More generally, we can replace $f(u)$ by a function $\hat{f}$ depending on $u,Du$ and satisfying a linear growth in $Du$. Namely, we will assume the following.

\bdes \item[F')] There exist a constant $C$ and a function $f(u)$ satisfying F) such that \beqno{FUDU11}|\hat{f}(u,Du)| \le C\llg^\frac12(u)|Du| + f(u),\eeq \beqno{fUDU21}|f_u(u| \le C\llg(u).\eeq
\edes

Let us discuss some applications where the conditions A) and F) can be verified. In many models, $A(u), f(u)$ have polynomial growths in $u$. That is, there are nonnegative numbers $k,K$ such that $\llg(u)\sim (1+|u|)^k$ and $|f(u)|\sim (1+|u|)^K$. Here and in the sequel, we will write $a\sim b$ if there are two generic positive constants $C_1,C_2$ such that $C_1b \le a \le C_2b$.

Obviously, \mref{Fghyp} in A) holds for any $k\ge0$. Concerning F), we can take $F(u)=|u|^\frac{k+2}{2}$ and $P(u)=|u|^{k+1}$. It is clear that \mref{dfuu} and \mref{Pu} hold for such choice of $F,P$.  We also have $|f(u)||u|\le (1+|u|)^{K+1}\sim (F(u)+1)^2$ if $K=k+1$. Thus, \mref{fuu} is satisfied with $\eg_0$ being the coefficient of the highest power of $u$ in $f(u)$.

On the other hand, if $K\le k+2$ then \mref{fuuu} is satisfied because $$\frac{\llg^\frac12(u)|f(u)|}{P(u)+1}\le C\frac{(1+|u|)^{\frac k2+K}}{(1+|u|)^{k+1}} \sim (1+|u|)^{K-\frac k2-1}\le C(1+|u|)^{\frac k2+1}=C(F(u)+1).$$ 

Hence, it is clear that the main assumptions in F) and \mref{fUDU21} in F') are verified if $K=k+1$. 

We then prove the following next main result which shows that if the {\bf excess}
\beqno{cca}\ccA(u)=\myprod{A(u)Du,DA(u)u_t-A(u)_tDu}\eeq of a classical solution $u$ to \mref{e1} does not blow up in time then \mref{blowup0} holds to give the global existence of the classical solution $u$.

\btheo{thm1} Assume A) and F'). Let $u$ be a classical solution to \mref{e1} and $T_0$ be its maximal existence time. 
For any $t_0\in(0,T_0)$ and $\eg>0$ assume that there are continuous function $C_\eg,C_\ccA$ on $(0,\infty)$ which may also depend on $\|u(\cdot,t_0)\|_{W^{1,2}(\Og)}$ and $\|u_t(\cdot,t_0)\|_{L^{2}(\Og)}$ such that  \beqno{AAu}\int_{t_0}^{s}\iidx{\Og}{\ccA(u)}\,dt \le \int_{t_0}^{s}\iidx{\Og}{[\eg\llg(u)|u_t|^2 +C_\eg(t)|A(u)Du|^2]}\,dt
+C_\ccA(s)\eeq for all $s\in(t_0,T_0)$.

If either $\eg_0$ or $d(\Og)$ is sufficiently small then $T_0=\infty$.

\etheo

The next two theorems show that \mref{AAu} holds for the generalized (SKT) and {\em pyramid} models.
\btheo{thm2} Assume A) and F'). Suppose that there is a $C^2$ function $\ccP:\RR^m\to\RR^m$ such that $A(u)=\ccP_u$.

If either $\eg_0$ or $d(\Og)$ is sufficiently small then $T_0=\infty$.

\etheo

\btheo{thm3}  For $k_0=1,\ldots,m-1$ we suppose that the subsystems of the first $k_0$ equations of \mref{e1} do not depend on the unknowns $u_i$ if $i>k_0$ and they satisfy A) and F'). Furthermore, there is a constant $C$ such that for $u=(u_1,\ldots,u_{k_0})^T$ \beqno{ak0}\left|\frac{\partial}{\partial u_{k_0}}a_{k_0j}(u)\right| \le C\llg^\frac12(u)\quad \forall j<k_0,\, k_0=2,\ldots,m.\eeq

If either $\eg_0$ or $d(\Og)$ is sufficiently small then $T_0=\infty$.

\etheo

We remmark that if $A(u)$ has polynomial growth in $u$ and $\llg(u)\sim (1+|u|)^k$ then \mref{ak0} holds if $k\in[0,2]$. In fact, the reaction terms of the subsystem for the first $k_0$ preys ($u_i$, $i\le k_0$) can depend on the predators $u_i$ ($i>k_0$). We are interested in the effect of cross diffusion in this paper and assume very weak feeding rate of the predators in the reaction terms of the systems. Strong feeding rates of predators for cross diffusion systems will be reported in our forthcoming works. 

\section{Proof}\eqnoset\label{proof}

In this section, we will consider a classical solution to \mref{e1} that exists in its maximal time interval $(0,T_0)$. We also fix a $t_0$ in $(0,T_0)$. We will frequently use functions $\ccC$ satisfying the following property.

\bdes \item[C)] $\ccC$ is continuous on $(0,\infty)$ and depends on $\|u(\cdot,t_0)\|_{W^{1,2}(\Og)}$, $\|(u(\cdot,t_0))_t\|_{L^{2}(\Og)}$ and $d(\Og)$.\edes

In the proof,  when there is no ambiguity $C, C_i$ will denote universal constants that can change from line to line in our argument. Furthermore, $C(\cdots)$ is used to denote quantities which are bounded in terms of theirs parameters. The same convention applies to functions $\ccC,\ccC_i$ satisfying the property C).

\blemm{ldulemm} Assume \mref{dfuu} and \mref{fuu} of F). If either $\eg_0$ or $d(\Og)$ is sufficiently small then there is a function $\ccC$ satisfying C) such that such that
\beqno{p00}\sup_{t\in[t_0,T]}\iidx{\Og}{|u|^2} + \itQ{\Og\times[t_0,T]}{\llg(u)|Du|^2} \le \ccC(T)\quad \forall T\in(t_0,T_0).\eeq
\elemm

\bproof Testing the system of $u$ with $u$ and integrating over $\Og\times[t_0,T]$ by parts, we easily obtain for $Q=\Og\times[t_0,T]$
\beqno{p0}\sup_{t\in[t_0,T]}\iidx{\Og}{|u|^2} + \itQ{Q}{\llg(u)|Du|^2} \le \itQ{Q}{\myprod{f(u),u}} +\|u(\cdot,t_0)\|_{L^2(\Og)}^2.\eeq

By \mref{fuu} of F) we can find a constant $C$ such that $$\iidx{\Og}{\myprod{f(u),u}} \le \eg_0\iidx{\Og}{F^2(u)} +C|\Og|.$$ Here, $|\Og|$ is the Lebesgue measure of $\Og$.

Because $F(0)=0$ and the boundary condition of $u$, we have $F(u)=0$  on the boundary $\partial\Og$. Using the Poincar\'e inequality and \mref{dfuu}, we have
$$\iidx{\Og}{F^2(u)} \le Cd^2(\Og)\iidx{\Og}{|DF(u)|^2} \le Cd^2(\Og)\iidx{\Og}{\llg(u)|Du|^2}.$$

Hence,
$$\iidx{\Og}{\myprod{f(u),u}} \le  C\eg_0d^2(\Og)\iidx{\Og}{\llg(u)|Du|^2} + C|\Og|.$$

Thus, if either $\eg_0$ or $d(\Og)$ is sufficiently small then \mref{p0} and the above yield a constant $C_1$ such that $$\sup_{t\in[t_0,T]}\iidx{\Og}{|u|^2} + C_1\itQ{Q}{\llg(u)|Du|^2} \le C|\Og|T+\|u(\cdot,t_0)\|_{L^2(\Og)}.$$ The last quantity defines a function $\ccC(T)$ satisfying C). The proof is then complete. \eproof

To proceed, we need the following elementary fact. By the Ladyzhenskyas inequality ($n=2$)  we have \beqno{lady}\left(\iidx{\Og}{|U|^4}\right)^\frac12\le \left(\iidx{\Og}{|U|^2}\right)^\frac12\left(\iidx{\Og}{|DU|^2}\right)^\frac12,\eeq
if $U=0$ on the boundary $\partial\Og$. The Poincar\'e inequality applies to the first factor on the right gives $$\left(\iidx{\Og}{|U|^4}\right)^\frac12  \le \iidx{\Og}{|DU|^2}.$$

Hence, if $U,V$ vanish on the boundary $\partial\Og$ then
\beqno{LP}\iidx{\Og}{|U|^2|V|^2} \le \left(\iidx{\Og}{|U|^4}\right)^\frac12
\left(\iidx{\Og}{|V|^4}\right)^\frac12\le C\iidx{\Og}{|DU|^2}\iidx{\Og}{|DV|^2}.\eeq

We also note that $\llg(u)$ is the smallest eigenvalue of $(A+A^T)/2$ and $\LLg(u)$ is the smallest eigenvalue of $A^TA$. Thus, if $\mu(u)$ is the eigenvalue of $A$ with smallest real part then $\llg(u)=\Re(\mu(u))$ and $\LLg(u)=|\mu(u)|^2$. Therefore,
\beqno{A2}|A(u)\zeta|^2=\myprod{A^T(u)A(u)\zeta,\zeta}\ge \LLg(u)|\zeta|^2 \ge \llg^2(u)|\zeta|^2.\eeq

\bproof {\em (Proof of \reftheo{thm1})} For any $t_0$ in $(0,T_0)$, our main goal is to show that there is a function $\ccC$ satisfying C) such that  \beqno{pduest} \iidx{\Og}{|A(u)Du(x,t)|^2} \le \ccC(t) \quad \forall t\in(t_0,T_0).\eeq 

Once this is established, by \mref{A2} and the fact that $\llg(u)$ is bounded from below by $\llg_0>0$, the integral of $\llg(u)|Du|^2$ over $\Og$ does not blow up in finite time. The condition \mref{blowup0} of \reftheo{thm0} follows and our theorem is then proved.

For any $T\in (t_0,T_0)$, test the system for $u$ by $A(u) u_t$ (i.e. multiplying the $i^{th}$ equation of \mref{e1} by $\sum_j a_{ij}(u_j)_t$, integrating over $\Og\times[t_0,T]$, summing the results) and integrate by parts to get
\beqno{p1}\itQ{\Og\times(t_0,T)}{(\myprod{A(u)u_t,u_t} +\myprod{A(u)Du,D(A(u)u_t)})}=\itQ{\Og\times(t_0,T)}{\myprod{f(u),A(u)u_t}}.\eeq

We note that $$\myprod{A(u)Du,D(A(u)u_t)} = \myprod{A(u)Du,DAu_t} +\myprod{A(u)Du,AD(u_t)},$$
$$\frac12\frac{\partial }{\partial t}\|ADu\|^2 = \myprod{A(u)Du,(A(u)Du)_t)}=\myprod{A(u)Du,A(u)_tDu}+\myprod{A(u)Du,ADu_t)}.$$

Hence, by the definition \mref{cca} of $\ccA$ $$\myprod{A(u)Du,D(A(u)u_t)}-\frac12\frac{\partial }{\partial t}|ADu|^2 =\myprod{A(u)Du,DAu_t-A(u)_tDu}=\ccA(u).$$

Thus,  we rewrite \mref{p1} as
$$\itQ{\Og\times(t_0,T)}{(\myprod{A(u)u_t,u_t} +\frac12\frac{\partial }{\partial t}|ADu|^2)}=\itQ{\Og\times(t_0,T)}{(\myprod{f(u),A(u)u_t}+\ccA(u))}.$$

The ellipticity of $A(u)$ and integrating in $t$ then give for any $T\in(t_0,T_0)$
$$\barr{ll}\lefteqn{\itQ{\Og\times(t_0,T)}{\llg(u)|u_t|^2}+ \frac12\iidx{\Og}{|A(u(x,T))Du(x,T)|^2}\le}\hspace{2cm}&\\& \frac12\iidx{\Og}{|A(u(x,t_0))Du(x,t_0)|^2}+\itQ{\Og\times(t_0,T)}{(C|f(u)|\llg(u)|u_t|+\ccA(u))}.\earr
$$

Using Young's inequality to find a constant $C(\eg)$ such that for any $\eg>0$
$$|f(u)|\llg(u)|u_t| \le \eg\llg(u)|u_t|^2 + C(\eg)\llg(u)|f(u)|^2.$$

For sufficiently small and fixed $\eg$ we then have 
\beqno{p2}\barr{ll}\lefteqn{\itQ{\Og\times(t_0,T)}{\llg(u)|u_t|^2}+  \iidx{\Og}{|A(u(x,T))Du(x,T)|^2}\le}\hspace{1cm}&\\& \iidx{\Og}{|A(u(x,t_0))Du(x,t_0)|^2}+C\itQ{\Og\times(t_0,T)}{(\llg(u)|f(u)|^2+\ccA(u))}.\earr
\eeq

Now, let $U=P(u)$ be the function in A) and $V=\frac{\llg^\frac12(u)|f(u)|}{P(u)+1}$. Then \mref{fuuu} gives $|V|\le CF(u)$ ($F(u)$ was also defined in A)). We observe that
$$\barr{ll}\lefteqn{\iidx{\Og}{\llg(u)|f(u)|^2}\le C\iidx{\Og}{(U^2+1)(F(u)^2+1)}}\hspace{1cm}&\\&= C\iidx{\Og}{P(u)^2F(u)^2} +C\iidx{\Og}{P(u)^2}+C\iidx{\Og}{F(u)^2} + C(d(\Og))\\ &\le C\iidx{\Og}{|DP(u)|^2}\iidx{\Og}{|DF(u)|^2}+C\iidx{\Og}{(|DP(u)|^2+|DF(u)|^2)} + C(d(\Og)),\earr$$ where we used \mref{LP} and then Poincar\'e's inequality for $P(u), F(u)$ in the last estimate.  By \mref{Pu} and \mref{A2}, we have  $$|DP(u)|^2 \le |P_u(u)|^2|Du|^2 \le C\llg^2(u)|Du|^2\le C\myprod{A^T(u)A(u)Du,Du} = C|A(u)Du|^2.$$

Since $|DF(u)|^2\le C\llg(u)|Du|^2$ by \mref{dfuu}, we can use the above estimates in \mref{p2} to get
\beqno{p3a}\barr{ll}\lefteqn{\itQ{\Og\times(t_0,T)}{\llg(u)|u_t|^2}+ \iidx{\Og}{|AD(u(x,T))|^2}\le}\hspace{1cm}&\\& \iidx{\Og}{|AD(u(x,t_0))|^2}+C\dspl{\int_{t_0}^{T}}\left[\iidx{\Og}{\llg(u)|Du|^2}+1\right]\iidx{\Og}{|A(u)Du|^2}\, dt\\&+C\dspl{\int_{t_0}^{T}}\iidx{\Og}{\llg(u)|Du|^2}\, dt+C(d(\Og))(T-t_0)+C\dspl{\int_{t_0}^{T}}\iidx{\Og}{\ccA(u)}\, dt.\earr
\eeq

By \mref{p00} and the assumption \mref{AAu} on $\ccA(u)$, There is a continuous function $C_\ccA(t)$ satisfying C) such that
$$\dspl{\int_{t_0}^{T}}\iidx{\Og}{\ccA(u)}\, dt \le\int_{t_0}^{T}\iidx{\Og}{[\eg\llg(u)|u_t|^2 +C_\eg(t)|A(u)Du|^2]}\,dt
+C_\ccA(T).$$

For a sufficiently small and fixed $\eg$ we derive from the above and \mref{p3a} that
\beqno{p3}\barr{ll}\lefteqn{\iidx{\Og}{|DA(u(x,T))|^2}\le \iidx{\Og}{|DA(u(x,t_0))|^2}+C_\ccA(T)}\hspace{2cm}&\\& +C\dspl{\int_{t_0}^{T}}[\iidx{\Og}{(\llg(u)|Du|^2+1)}+C_\eg(t)]\iidx{\Og}{|AD(u)|^2}\, dt.\earr
\eeq

We now set $$y(t)=\iidx{\Og}{|A(u)Du(x,t)|^2},\, \ag(t)=\iidx{\Og}{|A(u)Du(x,t_0)|^2}+C_\ccA(t),$$ and $$ \bg(t)=\iidx{\Og}{(\llg(u)|Du(x,t)|^2+1)} + C_\eg(t).$$

We obtain from \mref{p3} $$y(t) \le \ag(t)+C\int_{t_0}^{t}\bg(s)y(s)ds \quad \forall t\in(t_0,T_0).$$

The integral form of Gronwall's inequality gives
$$y(t) \le \ag(t)+C\int_{t_0}^{t}\ag(s)\bg(s)\exp\left(\int_{s}^{t}\bg(\tau)d\tau\right)ds.$$

Cleraly, there are functions $\ccC_1,\ccC_2$ satisfying C) such that $\ag(t)$ is bounded by $\ccC_1(t)$ and by \reflemm{ldulemm} $$\int_{s}^{t}\bg(\tau)d\tau \le  \itQ{\Og\times[t_0,t]}{\llg(u)|Du|^2} + Ct+\itQ{\Og\times[t_0,t]}{C_\eg(t)}\le \ccC_2(t) \quad \forall t\in(t_0,T_0).$$

We conclude that there is a function $\ccC_3$ satisfying C) such that $$y(t)=\iidx{\Og}{|A(u)Du(x,t)|^2} \le \ccC_3(t) \quad \forall t\in(t_0,T_0).$$ This gives the desired estimate \mref{pduest} for $Du$. The proof is complete. \eproof

\brem{fdurem} If we asume F') and replace $f(u)$ by $\hat{f}(u,Du)$ satisfying $$|\hat{f}(u,Du)| \le C\llg^\frac12(u)|Du| + f(u)$$ then the result continue to hold. Firstly, by Young's inequality $$\myprod{\hat{f}(u,Du),u} \le \eg \llg(u)|Du|^2 + C(\eg)|u|^2 + |f(u)||u|.$$  For suficiently small $\eg$, the argument in the proof of \reflemm{ldulemm} will lead to (the last inequality in the proof with an extra term)
$$\sup_{t\in[t_0,T]}\iidx{\Og}{|u|^2} \le  \ccC(T) + C(\eg)\itQ{\Og\times[0,T]}{|u|^2},$$ for some function $\ccC$ satisfying C). A simple use of Gronwall's inequality gives $$\sup_{t\in[t_0,T]}\iidx{\Og}{|u|^2} \le  C(T,\|u(\cdot,t_0)\|_{L^2(\Og)},d(\Og)),$$ and the assertion of the lemma still holds.

Next, $$\barr{lll}|\hat{f}(u,Du)|\llg(u)|u_t| &\le& C\llg(u)\llg(u)|Du||u_t| + C|f(u)|\llg(u)|u_t|\\&\le&\eg\llg(u)|u_t|^2 + C(\eg)\llg^2(u)|Du|^2 + C|f(u)|\llg(u)|u_t|.\earr$$

As $f(u)$ satisfies F), for small $\eg$ in the above the proof of \reftheo{thm1} can continue.

\erem 

We now consider the excess $\ccA$. If $A(u)=[a_{ij}(u)]$ then calculations give
$$A(u)Du = \left[\sum_{j}a_{ij}(u)Du_j\right],$$
$$DA(u)u_t = \left[\sum_k\frac{\partial}{\partial u_k}a_{ij}(u)Du_k\right]u_t=\left[\sum_{k,j}\frac{\partial}{\partial u_k}a_{ij}(u)Du_k(u_j)_t\right],$$
$$A(u)_tDu = \left[\sum_k\frac{\partial}{\partial u_k}a_{ij}(u)(u_k)_t\right]Du=\left[\sum_{k,j}\frac{\partial}{\partial u_k}a_{ij}(u)(u_k)_tDu_j\right].$$

Thus, $$DA(u)u_t-A(u)_tDu = \left[\sum_{k,j}\frac{\partial}{\partial u_k}a_{ij}(u)[Du_k(u_j)_t-(u_k)_tDu_j]\right].$$

In general, we have $$\barr{lll}\ccA(u)&=&\myprod{A(u)Du,DA(u)u_t-A(u)_tDu}\\&=&\sum_i\sum_{l}a_{il}(u)Du_l\sum_{k,j}\frac{\partial}{\partial u_k}a_{ij}(u)[Du_k(u_j)_t-(u_k)_tDu_j]\\&=&\sum_i\sum_{l,k,j}a_{il}(u)\frac{\partial}{\partial u_k}a_{ij}(u)[Du_k(u_j)_t-(u_k)_tDu_j]Du_l.\earr$$

Applying \reftheo{thm1} to special structures of $A(u)$ we have the proof of the last two theorems.

\bproof {\em (Proof of \reftheo{thm2})} If $\frac{\partial}{\partial u_k}a_{ij}(u)=\frac{\partial}{\partial u_i}a_{ik}(u)$ for any $i$ and $k\ne j$ then it is clear that $\ccA(u)=0$. This is the case if there are functions $\ccP_i$'s such that $a_{ij}(u)=\partial_{u_j}\ccP_i(u)$. The condition \mref{AAu} is satisfied and our theorem follows from \reftheo{thm1}. \eproof

\bproof {\em (Proof of \reftheo{thm3})} 
We assume that the $i$-th equation of the system \mref{e1} does not depend on the unknowns $u_k$ if $k>i$. This means $a_{ik}=0$ and $\frac{\partial}{\partial u_k}a_{ij}(u)=\frac{\partial}{\partial u_k}f_{i}(u)=0$ if $k>i$. 

Therefore, for the subsystem of the first $k_0$ equations ($k_0\ge1$), we have, writing $u=(u_1,\ldots,u_{k_0})^T$

$$A_{k_0}(u)Du = \left[\sum_{j\le k_0}a_{ij}(u)Du_j\right],$$
$$DA_{k_0}(u)u_t = \left[\sum_{j,k\le k_0}\frac{\partial}{\partial u_k}a_{ij}(u)Du_k\right]u_t=\left[\sum_{k,j\le k_0}\frac{\partial}{\partial u_k}a_{ij}(u)Du_k(u_j)_t\right],$$
$$A_{k_0}(u)_tDu = \left[\sum_{k\le k_0}\frac{\partial}{\partial u_k}a_{ij}(u)(u_k)_t\right]Du=\left[\sum_{k,j\le k_0}\frac{\partial}{\partial u_k}a_{ij}(u)(u_k)_tDu_j\right].$$

Accordingly, the excess for the subsystem is
 $$\barr{lll}\ccA_{k_0}(u)&=&\sum_{i,l,k,j\le k_0}a_{il}(u)\frac{\partial}{\partial u_k}a_{ij}(u)[Du_k(u_j)_t-(u_k)_tDu_j]Du_l\\
 &=&I_{1,k_0}+I_{2,k_0}+I_{3,k_0}.\earr$$

Here,
$$I_{1,k_0}=\sum_{j<k_0}a_{k_0k_0}(u)\frac{\partial}{\partial u_{k_0}}a_{k_0j}(u)[Du_{k_0}(u_j)_t-(u_{k_0})_tDu_j]Du_{k_0},$$
$$I_{2,k_0}=\sum_{j,l<k_0}a_{k_0l}(u)\frac{\partial}{\partial u_{k_0}}a_{k_0j}(u)[Du_{k_0}(u_j)_t-(u_{k_0})_tDu_j]Du_l,$$ and $$I_{3,k_0}=\sum_{i,j,k,l< k_0}a_{il}(u)\frac{\partial}{\partial u_k}a_{ij}(u)[Du_k(u_j)_t-(u_k)_tDu_j]Du_l.$$

 It is clear that $\ccA_1=0$ so that $u_1$ exists globally and its derivatives do not blow up in finite time. We then argue by induction. Consider the induction hypothesis:
 
 For some $k_0\ge1$ there is a continuous function $\ccC$ satisfying C) such that
 \beqno{indk0} \mbox{{\bf (K0)}} \quad \|(u_i(\cdot,t))_t\|_{L^\infty(\Og)}, \|u_i(\cdot,t)\|_{C^1(\Og)} \le \ccC_{k_0}(t)\quad \forall t\in (0,T_0),\, i=1,\ldots, k_0-1.\eeq

As we assume there is a constant $C$ such that $|A(u)|\le C\llg(u)$ and (see \mref{ak0})  $$\left|\frac{\partial}{\partial u_{k_0}}a_{k_0j}(u)\right| \le C\llg^\frac12(u)\quad \forall j<k_0,\, k_0=2,\ldots,m.$$

Then by Young's inequality and our induction assumptiom \mref{indk0} $$I_{1,k_0},I_{2,k_0} \le \eg\llg(u)|u_t|^2 + C(\eg)\ccC(t)\llg^2(u)|Du|^2 +C(\eg)\ccC(t)  \quad \forall t\in(t_0,T_0)$$ for some continuous function $\ccC$ depending on $\ccC_{k_0}$. Similarly, since $I_{3,k_0}$ does not depend on $u_{k_0}$, \mref{indk0} also shows  that $$\itQ{\Og\times(t_0,t)}{I_{3,k_0}} \le \ccC(t) \quad \forall t\in(t_0,T_0).$$

By \mref{A2}, $\llg^2(u)|Du|^2\le |A(u)Du|^2$. Therefore, we can conclude from the above estimates that there are continuous functions $\ccC, C_\ccA$ satisfying C) such that $$\int_{t_0}^{s}\iidx{\Og}{\ccA(u)}\,dt \le \int_{t_0}^{s}\iidx{\Og}{[\eg\llg(u)|u_t|^2 +C(\eg)\ccC(t)|A(u)Du|^2]}\,dt
+C_\ccA(s)$$ for all $s\in(t_0,T_0)$. Hence, \reftheo{thm1} can be applied to the subsystem of the first $k_0$ equations and we see that $u=(u_1,\ldots,u_{k_0})^T$ does not blow up in finite time. Therefore, \mref{indk0} holds again for $k_0+1$ and our proof is complete by induction. \eproof

\section{Proof of \reftheo{thm0}}\eqnoset\label{secthm0}

In the sequel, we will denote $\Fg(u)=\frac{|\llg_u(u)|^2}{\llg(u)}$. Before going to the proof, we need some estimates for the integral of $\Fg(u)|Du|^{4p}$.
 \blemm{Du4plemm} Assume \mref{Fghyp} in A). For any $p\ge1$ and any nonnegative function $\psi\in C^1_0(B_R)$ there is a constant $C_1$ such that
 \beqno{keyp}\barr{ll}\lefteqn{\iidx{B_R}{\Fg(u)|Du|^{4p}\psi^4}\le C\iidx{\Og}{\llg^2(u)|Du|^{4p}\psi^4} \le C_1\iidx{B_R}{\llg(u)|Du|^{2p}\psi^2}\times}\hspace{.5cm}&\\& \iidx{B_R}{(\llg(u)|Du|^{2p-2}|D^2u|^2\psi^2+\Fg(u)|Du|^{2p+2}\psi^2 + \llg(u)|D\psi|^2|Du|^{2p})}.\earr\eeq
  \elemm

 \bproof 
  By Ladyzhenskaya's inequality \mref{lady} with $U=\llg^\frac12(u)\psi|Du|^{p-1}Du$ we have $$\barr{lll}\iidx{B_R}{\llg^2(u)|Du|^{4p}\psi^4}&=&\iidx{B_R}{|U|^4} \le C\iidx{B_R}{|U|^2}\iidx{B_R}{|DU|^2}\\&\le&C\iidx{B_R}{\llg(u)|Du|^{2p}\psi^2}\iidx{B_R}{|D(\llg^\frac12(u)\psi|Du|^{p-1}Du)|^2}.\earr$$
 
 It is clear that there is a constant $C_2$ such that $$\barr{ll}\lefteqn{|D(\llg^\frac12(u)\psi|Du|^{p-1}Du)|^2 \le}\hspace{2cm}&\\& C_2\left[\llg(u)|Du|^{2p-2}|D^2u|^2\psi^2 + \frac{|\llg_u(u)|^2}{\llg(u)}|Du|^{2p+2}\psi^2 + \llg(u)|D\psi|^2|Du|^{2p}\right].\earr$$ 
 
Because $\llg(u)$ is bounded from below, \mref{Fghyp} gives $\Fg(u)\le C\llg^2(u)$ for some constant $C$. Applyling the above inequality in the previous estimate, we obtain the lemma. \eproof

Since $u$ are $C^2$ in $x$, we can differentiate \mref{e1} with respect to $x$ to get
\beqno{ga2} (Du)_t=\Div((A(u)D^2u + A_u(u)DuDu)+Df(u,Du).\eeq

Furthermore, by \cite[Lemma 2.1]{sd}, if $A$ is a matrix satisfying $\llg_0|\zeta|^2 \le \myprod{A\zeta,\zeta}$ and $|A\zeta|\le \LLg_0|\zeta|$ then for any $\ag$ and $\dg_\ag\in(0,1)$ are numbers such that  $\frac{\ag}{2+\ag}=\dg_{\ag}\frac{\llg_0}{\LLg_0}$ then there is a positive constant $\hat{\llg}$ depending on $\llg_0,\LLg_0,\dg_\ag$ such that
 \beqno{D2v}\myprod{AD\zeta, D(\zeta |\zeta |^\ag)}\ge
 \widehat{\llg}|\zeta |^{\ag}|D\zeta |^2.\eeq
 
 We also recall the following elementary iteration result (e.g., see \cite[Lemma 6.1, p.192]{Gius}). 
 
 \blemm{Giusiter} Let $f,g,h$ be bounded nonnegative functions in the interval $[\rg,R]$ with $g,h$ being increasing. Assume that for $\rg \le s<t\le R$ we have $$f(s) \le [(t-s)^{-\ag} g(t)+h(t)]+\eg f(t)$$ with $C\ge0$, $\ag>0$ and $0\le\eg<1$. Then $$f(\rg) \le c(\ag,\eg)[(R-\rg)^{-\ag} g(R)+h(R)].$$ The constant $c(\ag,\eg)$ can be taken to be $(1-\nu)^{-\ag}(1-\nu^{-\ag}\nu_0)^{-1}$ for any $\nu$ satisfying $\nu\in(0,1)$ and $\nu^{-\ag}\nu_0<1$.\elemm

\bproof  Let $B_R=B_R(x_0)$ be a ball in $\Og$. We consider only the interior case because the boundary case, when the center $x_0$ of $B_R$ is on the boundary $\partial\Og$, is similar.
For any $s,t$ such that $0\le s<t \le R$ let $\psi$ be a cutoff function for two balls $B_s,B_t$ centered at $x_0$. That is, $\psi\equiv1$ in $B_s$ and $\psi\equiv0$ outside $B_t$ with $|D\psi|\le1/(t-s)$. In the sequel, we will fix two reals $t_0,T$ such that $0<t_0<T<T_0$.

We divide the proof in three steps.

{\bf Step 1:} (Local energy estimates)  Clearly, by uniform ellipticity of $A(u)$, we can find a constant $C_0$ such that $|A(u)\zeta|\le C_0\llg(u)|\zeta|$. Thus, for any $p>1$ there is $\dg_p\in(0,1)$ such that $\ag=2p-2$ satisfies \beqno{pcond} \frac{\ag}{2+\ag}=\frac{2p-2}{2p}=\dg_{p}C_0^{-1}=\dg_{p}\frac{\llg(u)}{C_0\llg(u)}.\eeq
 
Testing \mref{ga2} with $|Du|^{2p-2}Du\psi^2$. 
By the above, there is a positive constant $C(p)$ (see \mref{D2v}) such that for $Q=\Og\times[t_0,T_0]$
$$\barr{ll}\lefteqn{\sup_{\tau\in(t_0,T)}\iidx{\Og}{|Du|^{2p}\psi^2}+C(p)\itQ{Q}{\llg(u)|Du|^{2p-2}|D^2u|^2\psi^2} \le}\hspace{.2cm}&\\& \itQ{Q}{|A(u)||D^2u||Du|^{2p-1}\psi|D\psi|}-\itQ{Q}{A_u(u)DuDuD(|Du|^{2p-2}Du\psi^2)}\\& + \itQ{Q}{D\hat{f}(u,Du)|Du|^{2p-2}Du\psi^2}+\iidx{\Og}{|Du(x,t_0)|^{2p}\psi^2}.\earr$$

For simplicity, we will assume in the sequel that $\hat{f}\equiv0$. The presence of $\hat{f}$ will be discussed later in \refrem{fdurem1}. For any given positive $\eg$ we use Young's inequality to find a constant $C(\eg)$ such that
$$|A(u)||D^2u||Du|^{2p-1}\psi|D\psi| \le \eg \llg(u)|Du|^{2p-2}|D^2u|^2\psi^2 + C(\eg)\llg(u)|Du|^{2p}|D\psi|^2,$$
$$\barr{ll}\lefteqn{|A_u(u)DuDuD(|Du|^{2p-2}Du\psi^2)|\le |A_u(u)||Du|^{2p}|D^2u|\psi^2+|A_u(u)||Du|^{2p+1}\psi|D\psi|}\hspace{2cm}&\\&\le \eg\llg(u)|Du|^{2p-2}|D^2u|^2 + C(\eg)\frac{|A_u|^2}{\llg(u)}|Du|^{2p+2}\psi^2 +C(\eg)\llg(u)|Du|^{2p}|D\psi|^2.\earr$$

Therefore, taking $\eg$ small and the above two inequalities in the previous one, we easily deduce  
\beqno{keypp}\barr{ll}\lefteqn{\sup_{\tau\in(t_0,T)}\iidx{B_s}{|Du|^{2p}}+\itQ{Q_s}{\llg(u)|Du|^{2p-2}|D^2u|^2} \le }\hspace{1cm}&\\&C_1\itQ{Q_t}{\Fg(u)|Du|^{2p+2}\psi^2}+\frac{C_1}{(t-s)^2}\itQ{Q_t}{\llg(u)|Du|^{2p}}+C(t_0).\earr\eeq  Here, we  used the defintion of $\psi$ and $\Fg(u)$ and the fact that $|A_u|\sim\llg_u$, and denoted $$Q_t=B_t\times(t_0,T),\, C(t_0)=\iidx{\Og}{|Du(x,t_0)|^{2p}}.$$

\newc{\ccG}{{\cal G}}

We now set $$\ccA_p(t)=\sup_{\tau\in(t_0,T)}\iidx{B_t}{|Du|^{2p}},\, \ccH_p(t) = \itQ{Q_t}{\llg(u)|Du|^{2p-2}|D^2u|^2},$$
$$\ccB_p(t)= \itQ{Q_t}{\Fg(u)||Du|^{2p+2}},\,\ccG_p(t)= \itQ{Q_t}{\llg(u)||Du|^{2p}} +C(t_0).$$

So that, for any $p\ge1$ and satisfies \mref{pcond} and $t$ is small, \mref{keypp} can be rewritten as
\beqno{keyppp} \ccA_p(s) + \ccH_p(s) \le C_1\ccB_p(t) + \frac{C_1}{(t-s)^2}\ccG_p(t) \quad \forall s,t \mbox{ such that } 0<s<t\mbox{ and } B_t\subset\Og.\eeq

{\bf Step 2:} (Estimates for the integral of $|Du|^4$ over $\Og\times(0,T)$) By \mref{blowup0} and the uniform continuity of integrals (see \refrem{hsv} following the proof) give for any given $\eg_0>0$ a constant $R(\eg_0,T)$ such that \beqno{hsvT}\iidx{B_R}{\llg(u)|Du(x,\tau)|^2} \le \eg_0 \quad \forall R< R(\eg_0,T),\, \tau\in[t_0,T].\eeq Therefore, let $p=1$ and $\psi$ be the cutoff function for $B_s,B_t$  in \mref{keyp} and use the definition of $\psi,\ccB_p,\ccH_p$ and $\ccC_p$ to have by integrating in $[t_0,T]$
$$\ccB_1(s) \le \sup_{\tau\in[t_0,T]}\iidx{B_R}{\llg(u)|Du(x,\tau)|^2}\left(\ccH_1(t)+\ccB_1(t) + \frac{1}{(t-s)^2}\ccG_1(t)\right)$$ for all $s,t$ such that $0<s<t<R(\eg_0,T)$.
We now choose $\eg_0$ suficiently small in \mref{hsvT} to have a number $\mu_0\in(0,1)$ such that
\beqno{x1}C_1\ccB_1(s) \le \frac{\mu_0}{2}\left(\ccH_1(t)+\ccB_1(t) + \frac{1}{(t-s)^2}\ccG_1(t)\right).\eeq

For $p=1$, \mref{keyppp} gives $$\ccH_1(s) \le C_1 \ccB_1(t) + \frac{C_1}{(t-s)^2}\ccG_1(t),\quad 0<s<t<R(\eg_0,T).$$ Let $t_1=(s+t)/2$ and use \mref{x1} with $s$ being $t_1$ and the above with $t$ being $t_1$ to obtain
\beqno{x2} \ccH_1(s) \le \frac{\mu_0}{2}[\ccH_1(t)+\ccB_1(t)] + \frac{C_2}{(t-s)^2}\ccG_1(t).\eeq

Obviously, we can assume $C_1\ge1$ so that we can add \mref{x1} and \mref{x2} to have $$\ccH_1(s) +\ccB_1(s) \le \mu_0 [\ccH_1(t)+\ccB_1(t)] + \frac{C_3}{(t-s)^2}\ccG_1(t),\quad 0<s<t<R(\eg_0,T).$$

Since $\mu_0\in(0,1)$, we can use \reflemm{Giusiter} with $f(t)=\ccH_1(t)+\ccB_1(t)$, $h(t)=0$, $g(t)=\ccG_1(t)$ and $\ag=2$ to obtain a constant $C_4$ depending on $\mu_0,C_3$ such that
$$\ccH_1(s)+\ccB_1(s) \le \frac{C_4}{(t-s)^2}\ccG_1(t),\quad 0<s<t<R(\eg_0,T).$$

For $R=R(\eg_0,T)/4$, the above with $s=R$, $t=2R$ gives \beqno{H1}\ccH_1(R)+\ccB_1(R) \le \frac{C_4}{R^2}\itQ{Q_{2R}}{\llg(u)|Du|^{2}}+C(t_0).\eeq

 Hence, by \mref{x1} and the estimate for the integral of $\llg(u)|Du|^2$ over $Q$,  we have \beqno{keydu44} \itQ{Q_R}{\Fg(u)|Du|^4} \le C(T,R(\eg_0),\|u(\cdot,t_0)\|_{C^1(\Og)})\quad \forall R< R(\eg_0,T).\eeq 
 We also note that \mref{H1}, \mref{keydu44} and the second inequality of \mref{keyp} gives
\beqno{keydu44z} \itQ{Q_R}{\llg^2(u)|Du|^4} \le C(T,R(\eg_0),\|u(\cdot,t_0)\|_{C^1(\Og)})\quad \forall R< R(\eg_0,T).\eeq
 
 Finite covering of $\Og$ by balls $B_{R(\eg_0,T)/2}$ yields \beqno{keydu444} \itQ{Q}{\llg^2(u)|Du|^4} \le C(T,R(\eg_0,T),\Og,\|u(\cdot,t_0)\|_{C^1(\Og)}) \quad \forall T\in(0,T_0).\eeq

{\bf Step 3:} (Estimates for the integral of $|Du|^{2p}$ over $\Og$) For any $p>1$ we have by H\"older's inequality \beqno{H11}\itQ{Q_t}{\Fg(u)|Du|^{2p+2}\psi^2} 
\le\left(\itQ{Q_t}{\Fg^2(u)|Du|^{4}}\right)^\frac12\left(\itQ{Q_t}{|Du|^{4p}\psi^4}\right)^\frac12.\eeq

Using Ladyzhenskaya's inequality \mref{lady} with $U=|Du|^{p-1}Du\psi$ and integrating the result over $(t_0,T)$, we have
$$\itQ{Q_t}{|Du|^{4p}\psi^4}\le C\sup_{\tau\in(t_0,T)}\iidx{B_t}{|Du|^{2p}\psi^2}\itQ{Q_t}{|DU|^2}$$

Since $\llg(u)$ is bounded from below by $\llg_0$, there is a constant $C(\llg_0)$ such that $$|DU|^2 \le C(\llg_0)[\llg(u)|Du|^{2p-2}|D^2u|^2\psi^2+\llg(u)|Du|^{2p}|D\psi|^2].$$
 Therefore,
$$\itQ{Q_t}{|Du|^{4p}\psi^4}\le C\ccA_p(t)[\ccH_p(t)+\frac{1}{(t-s)^2}\ccG_p(t)] \le C[\ccA_p(t) + \ccH_p(t)+\frac{1}{(t-s)^2}\ccG_p(t)]^2.$$

Here, Cauchy's inequality was used in the last inequality. Using the above estimate in \mref{H11} and the fact that $\Fg^2(u)\le C\llg^2(u)$, we derive
$$\itQ{Q_t}{\Fg(u)|Du|^{2p+2}\psi^2} \le C\left(\itQ{Q_t}{\llg^2(u)|Du|^{4}\psi^4}\right)^\frac12[\ccA_p(t) + \ccH_p(t)+\frac{1}{(t-s)^2}\ccG_p(t)].$$ 

By \mref{keydu444} and the continuity of integrals, the first factor on the right can be as small as we please if $t$ is small. Hence, for any given $\mu_1\in(0,1)$,  and the definition of $\ccB_p$, if $R\le R(\mu_1,T)$ for some small $R(\mu_1,T)$ then the above gives
$$C_1 \itQ{Q_t}{\Fg(u)|Du|^{2p+2}\psi^2} \le \mu_1[\ccA_p(t) + \ccH_p(t)+\frac{1}{(t-s)^2}\ccG_p(t)] \quad\mbox{ for some $\mu_1\in(0,1)$}.$$

If $p>1$ and satisfies \mref{pcond}, we then have from \mref{keypp} and the above inequality the following. $$\ccA_p(s) + \ccH_p(s) \le \mu_1 (\ccA_p(t)+\ccH_p(t)) + \frac{1}{(t-s)^2}\ccG_p(t), \quad 0<s<t<R(\mu_1,T).$$

For $f(t)=\ccA_p(t) + \ccH_p(t)$, $h(t)=0$, $g(t)=\ccG_p(t)$ and $\ag=2$ we can use \reflemm{Giusiter}, as $\mu_1\in(0,1)$,  to obtain
$$\ccF(\rg) \le\frac{C_5(\mu_1)}{(R-\rg)^2}\ccG_p(R), \quad 0<\rg<R <R(\mu_1,T).$$

We can assume that $2p<4$. Because $\llg(u)$ is bounded from below, by \mref{keydu44z} and a simple use of H\"older's inequality, we can see that $\ccG_p(R)$ is bounded, using \mref{keydu44z}. 
Hence, the above yields \beqno{Viter11} \sup_{t\in(t_0,T)}\iidx{B_\rg}{|Du|^{2p}} +\itQ{Q_\rg }{\llg(u)|Du|^{2p-2}|D^2u|^2}\le C(\rg,R,\|u(\cdot,t_0)\|_{C^1(\Og)})\eeq if $0<\rg<R$ and $R$ is sufficiently small and some $p>1$. Finite covering of $\Og$ with balls $B_{R/2} $ yields
\beqno{Viter1zz1} \sup_{t\in(t_0,T)}\iidx{\Og}{|Du|^{2p}}+\itQ{Q}{\llg(u)|Du|^{2p-2}|D^2u|^2} \le C(\Og,R,\|u(\cdot,t_0)\|_{C^1(\Og)}).\eeq

Since $2p>2$, Sobolev's imbedding theorem shows that $u$ is H\"older continuous in $x$. From the system for $u$ and the above, with $p=1$, we see that $u_t$ is in $L^2(Q)$. It is now standard to show that $u$ is H\"older in $(x,t)$ and $Du$ is H\"older continuous. We now can refer to Amann's results to see that $u$ exists globally.  \eproof

\brem{fdurem1} If we replace $f(u)$ by a function $\hat{f}$ depending on $u,Du$ and satisfying a linear growth in $Du$ then the proof can go on with minor modification. Namely, there exist a constant $C$ and a function $f(u)$ satisfying F) such that $$|\hat{f}(u,Du)| \le C\llg^\frac12(u)|Du| + f(u).$$ We can assume that $|D\hat{f}(u,Du)| \le C|D(\llg^\frac12(u)|Du|) + |f_u(u)||Du|$ so that $$|D\hat{f}(u,Du)| \le  C\llg^\frac12|D^2u|+ C\Fg^\frac12(u)|Du|^2+ |f_u(u)||Du|.$$ Therefore, in {\bf Step 1}, the extra term $|D\hat{f}(u,Du)||Du|^{2p-1}\psi^2$ can be handled by using the following estimates, which are the results of a simple use of Young's inequality.
$$\barr{lll}|D\hat{f}(u,Du)||Du|^{2p-1} &\le&  C[\llg^\frac12|D^2u|+ C\Fg^\frac12(u)|Du|^2+ |f_u(u)||Du|]|Du|^{2p-1}\\&\le&
\eg\llg|Du|^{2p-2}|D^2u|^2 + C(\eg)\llg|Du|^{2p} +\\
&&C\Fg(u)|Du|^{2p+2} + C|Du|^{2p} +C|f_u||Du|^{2p}\earr.$$

We can then assume that $|f_u|\le C\llg(u)$ for some constant $C$ and see that the proof can continue to obtain the energy estimate \mref{keypp}. The result then follows.

\brem{hsv} The existence of $R=R(\eg_0,T)$ in \mref{hsvT} is an easy consequence of  a simple application of Hahn-Saks-Vitali's theorem and an argument by contradiction. Indeed, if there is no such uniform $R$ for \mref{hsvT} to hold then there is a sequence $\{s_n\}$ in $[t_0,T]$ such that the integrals of $f_n:=\llg(u(\cdot,s_n))|Du(\cdot,s_n)|^2$ over $B_{\frac1n}\times\{t_n\}$ is greater than $\eg_0$. We can assume that $\{s_n\}$ converges in $[t_0,T]$. Since $\llg(u)|Du|^2$ is continuous in $t\in[t_0,T]$, Hahn-Saks-Vitali's theorem, e.g. see \cite{DS}, applies to the sequence $f_n$ and shows the uniform continuity in $R$ of the integral of $\llg(u)|Du|^2$ over $B_{R}\times\{t_n\}$ and gives a contradiction.
\erem

\erem 

\section{Further Discussion} \eqnoset\label{prevsec}
To compare \reftheo{thm0} with some results in our earlier work \cite{leGN}, where we considered general dimension $n\ge2$, let us recall the following results in \cite{leGN}. There, we assumed the following structural conditions on the system \mref{e1}

\bdes \item[A.1)] (Uniform ellipticity) There are positive constants $C,\llg_0$ and a smooth function $\llg(u)$ such that $\llg(u)\ge \llg_0$ and $$\llg(u)|\xi|^2 \le \myprod{A(u)\xi,\xi} \le C\llg(u)|\xi|^2\quad \forall u\in\RR^m,\,\xi\in\RR^{nm}.$$

\item[A.2)] Assume that $A\in C^1(\RR^m)$. Let $\Fg_0,\Fg$ be defined as $$\Fg_0(u) = \llg^\frac12(u)\mbox{ and } \Fg(u)= \frac{|A_u(u)|}{\llg^\frac12(u)}\quad u\in\RR^m.$$ Assume that the quantities
\beqno{CFg}
k_1:=\sup_{u\in\RR^m}\frac{|\Fg_u|}{\Fg},\, k_2:=\sup_{u\in\RR^m}\frac{\Fg}{\Fg_0}\eeq are finite.
\item[A.3)](Weights) If $u\in BMO(\Og)$ then $\Fg^\frac23(u)$ belongs to the $A_\frac43$ class and the quantity $[\Fg(u)^\frac23]_{\frac43}$ can be controlled by the norm $\|u\|_{BMO(\Og)}$.
\edes

As we discussed in \cite{leGN}, if $A(u)$ has a polynomial growth in $u$ then A.2) is easily satisfied. In the general case, our assumption \mref{Fghyp} in this paper is clearly much weaker than \mref{CFg}. Concerning the verification of A.3), a crucial factor for the validity of a generalized weighted Gagliardo-Nirenberg inequality involving BMO norms, extending a result in \cite{SR}, we used a connection between BMO functions and $A_\cg$ weights (see \cite{JN,OP}) to see that if $\llg(u)\sim (1+|u|)^k$ for $0\le k<5$ then A.3) holds. In this paper, we don't need such weighted Gagliardo-Nirenberg inequality so that we can allow $k$ to be any nonnegative number.

We also assumed in \cite{leGN} that the ellipticity constants $\llg(u),\LLg(u)$ of the matrix $A(u)$ were not too far apart.
\bdes
\item[R)] (The ratio condition) There is $\dg\in[0,1)$ such that \beqno{ratio} \frac{n-2}{n} = \dg \sup_{u\in\RR^m}\frac{\llg(u)}{\LLg(u)}.\eeq
\edes

One should note that there are examples in \cite{sd} of blow up solutions to \mref{e1} if the condition R) is violated. In this paper, when $n=2$, \mref{ratio} is clearly not needed.

We assumed in \cite{leGN} the following growth conditions on the nonlinearity $f$.

\bdes \item[F)] There are positive constants $C,b$ such that for any vector valued functions $u\in C^1(\Og,\RR^m)$ and $p\in C^1(\Og,\RR^{mn})$ \beqno{fcond2a} |f(u,p)| \le C|p| + C|u|^b+C.\eeq
\beqno{fcond3a} |Df(u,p)| \le C|Dp| + C|u|^{b-1}|Du|+C.\eeq
\edes

Under the above structural conditions, we proved the following global existence result of classical solutions and improved \reftheo{Amthm}. 

\btheo{Amthm1} (Theorem 2.5 in \cite{leGN})  Assume A.1)-A.3), R) and F).  Let $p_0 \in (n,\infty)$ and $U_0$ be  in $W^{1,p_0}(\Og
)$. Suppose that $T_0\in(0,\infty]$ is the maximal existence time for a classical solution $$u\in C([0,T_0),W^{1,p_0}(\Og)) \cap C^{1,2}((0, T_0)\times\bar{\Og})$$ for the system
\mref{e1}. 
Suppose that there is a function $C$ in $C^0((0,T_0])$ such that $$ \|u(\cdot,t)\|_{BMO(\Og)} \le C(t) \quad \forall t\in(0,T_0).$$
Moreover, for any $\eg>0$ and $(x,t)\in Q$, there exists $R=R(\eg)>0$ such that \beqno{vmouni1} \|u(\cdot,t)\|_{BMO(B_R(x))}<\eg\quad \forall t\in(0,T_0).\eeq

Then $T_0=\infty$.
\etheo

By Poincar\'e's inequality, the following consequence follows easily.

\bcoro{Amz} (Corollary 2.6 in \cite{leGN}) In addition to the assumptions of \reftheo{Amthm}, we assume R). Then there exists a maximal time $T_0\in(0,\infty]$ such that the system
\mref{e1} has a unique classical solution in $(0, T_0)$ with
$$u\in C([0,T_0),W^{1,p_0}(\Og)) \cap C^{1,2}((0, T_0)\times\bar{\Og})$$
Moreover, if $T_0<\infty$ then \beqno{blowup1}\lim_{t\to T^-_0}\|u(\cdot,t)\|_{W^{1,n}(\Og)}=\infty.\eeq\ecoro

As we discussed earlier, our structural conditions for \mref{e1} in this paper, when $n=2$, is much more general than those of \cite{leGN} to obtain the same conclusion of the above corollary in \reftheo{thm0}. Although \mref{blowup0} is slightly stronger than \mref{blowup1} but we have proved that it could be verified in many applications. In fact, it is possible to assume \mref{blowup1} and obtain \reftheo{thm0} under a bit stronger assumption than that of A) in this paper.

\bibliographystyle{plain}

\end{document}